\documentclass{article}
\usepackage{amsmath}

\newtheorem{theorem}{Theorem}

\newtheorem{example}[theorem]{Example}

\newtheorem{lemma}[theorem]{Lemma}

\newtheorem{remark}[theorem]{Remark}

\newenvironment{proof}[1][Proof]{\textbf{#1.} }{\ \rule{0.5em}{0.5em}}

\begin{document}

\title{Error inequalities for an optimal 3-point quadrature formula of
closed type }
\author{Nenad Ujevi\'{c}}
\maketitle

\begin{abstract}
An optimal 3-point quadrature formula of closed type is derived. Various
error inequalities are established. Applications in numerical integration
are also given.
\end{abstract}



\section{Introduction}

In recent years a number of authors have considered an error analysis for
quadrature rules of Newton-Cotes type. In particular, the mid-point,
trapezoid and Simpson rules have been investigated more recently (\cite{C1}, %
\cite{DAC1}, \cite{DCR1}, \cite{DPW1}, \cite{PPUV1}) with the view of
obtaining bounds on the quadrature rule in terms of a variety of norms
involving, at most, the first derivative. In the mentioned papers explicit
error bounds for the quadrature rules are given. These results are obtained
from an inequalities point of view. The authors use Peano type kernels for
obtaining a specific quadrature rule.

Quadrature formulas can be formed in many different ways. For example, we
can integrate a Lagrange interpolating polynomial of a function $f$ to
obtain a corresponding quadrature formula (Newton-Cotes formulas). We can
also seek a quadrature formula such that it is exact for polynomials of
maximal degree (Gauss formulas). Gauss-like quadrature formulas are
considered in \cite{U5}.

Here we present a new approach to this topic. Namely, we give a type of
quadrature formula. We also give a way of estimation of its error and all
parameters which appear in the estimation. Then we seek a quadrature formula
of the given type such that the estimation of its error is best possible.
Let us consider the above described procedure with more details.

If we define%
\begin{equation*}
K_{2}(\alpha ,\beta ,\gamma ,\delta ,t)=\left\{ 
\begin{array}{c}
\frac{1}{2}(t-\alpha )(t-\beta )\text{, }t\in \left[ a,\frac{a+b}{2}\right]
\\ 
\frac{1}{2}(t-\gamma )(t-\delta )\text{, }t\in \left( \frac{a+b}{2},b\right]%
\end{array}%
\right.
\end{equation*}%
then, integrating by parts, we obtain%
\begin{eqnarray*}
&&\int\limits_{a}^{b}K_{2}(\alpha ,\beta ,\gamma ,\delta ,t)f^{\prime \prime
}(t)dt \\
&=&\frac{1}{2}\left\{ f^{\prime }(b)\left( b-\gamma \right) \left( b-\delta
\right) +f^{\prime }\left( \text{$\frac{a+b}{2}$}\right) \times \right. \\
&&\left[ \left( \text{$\frac{a+b}{2}$}-\alpha \right) \left( \text{$\frac{a+b%
}{2}$}-\beta \right) -\left( \text{$\frac{a+b}{2}$}-\gamma \right) \left( 
\text{$\frac{a+b}{2}$}-\delta \right) \right] \\
&&\left. -f^{\prime }(a)(a-\alpha )(a-\beta )\right\} \\
&&+\left( a-\text{$\frac{\alpha +\beta }{2}$}\right) f(a)-\left( \text{$%
\frac{\gamma +\delta }{2}$-$\frac{\alpha +\beta }{2}$}\right) f\left( \text{$%
\frac{a+b}{2}$}\right) -\left( b-\text{$\frac{\gamma +\delta }{2}$}\right)
f(b) \\
&&+\int\limits_{a}^{b}f(t)dt.
\end{eqnarray*}

If we choose $\alpha =\beta =a$ and $\gamma =\delta =b$ then we get the
mid-point quadrature rule. If we choose $\alpha =\gamma =a$ and $\beta
=\delta =b$ then we get the trapezoid rule. If we choose $\alpha =0$, $\beta
=\frac{a+2b}{3}$ and $\gamma =\frac{2a+b}{3}$, $\delta =1$ then we get
Simpson's rule.

If we require that%
\begin{equation*}
\left( b-\gamma \right) \left( b-\delta \right) =0
\end{equation*}%
\begin{equation*}
\left( \text{$\frac{a+b}{2}$}-\alpha \right) \left( \text{$\frac{a+b}{2}$}%
-\beta \right) -\left( \text{$\frac{a+b}{2}$}-\gamma \right) \left( \text{$%
\frac{a+b}{2}$}-\delta \right) =0
\end{equation*}%
\begin{equation*}
(a-\alpha )(a-\beta )=0
\end{equation*}%
then we get a classical quadrature formula of the form%
\begin{eqnarray}
&&\int\limits_{a}^{b}K_{2}(\alpha ,\beta ,\gamma ,\delta ,t)f^{\prime \prime
}(t)dt  \label{j11as} \\
&=&\left( a-\text{$\frac{\alpha +\beta }{2}$}\right) f(a)-\left( \text{$%
\frac{\gamma +\delta }{2}-\frac{\alpha +\beta }{2}$}\right) f\left( \text{$%
\frac{a+b}{2}$}\right) -\left( b-\text{$\frac{\gamma +\delta }{2}$}\right)
f(b)  \notag \\
&&+\int\limits_{a}^{b}f(t)dt.  \notag
\end{eqnarray}

In practice we cannot find an exact value of the remainder term (error)

$\int\limits_{a}^{b}K_{2}(\alpha ,\beta ,\gamma ,\delta ,t)f^{\prime \prime
}(t)dt$. All we can do is to estimate the error. It can be done in different
ways. For example,

\begin{equation}
\left| \int\limits_{a}^{b}K_{2}(\alpha ,\beta ,\gamma ,\delta ,t)f^{\prime
\prime }(t)dt\right| \leq \underset{t\in \left[ a,b\right] }{\max }\left|
f^{\prime \prime }(t)\right| \int\limits_{a}^{b}\left| K_{2}(\alpha ,\beta
,\gamma ,\delta ,t)\right| dt.  \label{j12}
\end{equation}

It is a natural question which formula of the type (\ref{j11as}) is optimal,
with respect to a given way of estimation of the error. The main aim of this
paper is to give an answer to this question and to consider the formula from
an inequalities point of view. In fact, we seek a quadrature formula of the
given type such that its error bound is minimal. Note that we can minimize
only the factor $\int\limits_{a}^{b}\left| K_{2}(\alpha ,\beta ,\gamma
,\delta ,t)\right| dt$ in (\ref{j12}). A general approach is: we first
consider the minimization problem and then we formulate final results.
Vrious error inequalities for the obtained optimal formula are established.
Applications in numerical integration are also given. Finally, let us
mention that the obtained optimal quadrature formula has better estimations
of error than the Simpson's formula (see Remark \ref{Rs}).

\section{An optimal quadrature formula}

We consider the problem, described in Section 1, on the interval $\left[ 0,1%
\right] $. Let $\alpha ,\beta ,\gamma ,\delta \in R$. We define the mapping%
\begin{equation}
K_{2}(\alpha ,\beta ,\gamma ,\delta ,t)=\left\{ 
\begin{array}{c}
\frac{1}{2}(t-\alpha )(t-\beta )\text{,}t\in \left[ 0,\frac{1}{2}\right] \\ 
\frac{1}{2}(t-\gamma )(t-\delta )\text{, }t\in \left( \frac{1}{2},1\right]%
\end{array}%
\right. .  \label{j21n}
\end{equation}%
Let $I\subset R$ be an open interval such that $\left[ 0,1\right] \subset I$
and let $f:I\rightarrow R$ be a twice differentiable function such that $%
f^{\prime \prime }$ is bounded and integrable. We denote%
\begin{equation}
\left\| f\right\| _{\infty }=\underset{t\in \left[ 0,1\right] }{\sup }\left|
f(t)\right| .  \label{j22n}
\end{equation}%
Integrating by parts, we obtain%
\begin{eqnarray}
&&\int\limits_{0}^{1}K_{2}(\alpha ,\beta ,\gamma ,\delta ,t)f^{\prime \prime
}(t)dt  \label{j232a} \\
&=&\frac{1}{2}\int\limits_{0}^{\frac{1}{2}}(t-\alpha )(t-\beta )f^{\prime
\prime }(t)+\frac{1}{2}\int\limits_{\frac{1}{2}}^{1}(t-\gamma )(t-\delta
)f^{\prime \prime }(t)dt  \notag
\end{eqnarray}

\begin{eqnarray}
&=&-\frac{1}{2}\alpha \beta f^{\prime }(0)+\frac{1}{2}(1-\gamma )(1-\delta
)f^{\prime }(1)  \notag \\
&&+\frac{1}{2}\left[ \left( \frac{1}{2}-\alpha \right) \left( \frac{1}{2}%
-\beta \right) -\left( \frac{1}{2}-\gamma \right) \left( \frac{1}{2}-\delta
\right) \right] f^{\prime }(\frac{1}{2})  \notag \\
&&-\int\limits_{0}^{1}K_{1}(\alpha ,\beta ,\gamma ,\delta ,t)f^{\prime
}(t)dt,  \notag
\end{eqnarray}

where

\begin{equation*}
K_{1}(\alpha ,\beta ,\gamma ,\delta ,t)=\left\{ 
\begin{array}{c}
t-\frac{\alpha +\beta }{2},\text{ }t\in \left[ 0,\frac{1}{2}\right] \\ 
t-\frac{\gamma +\delta }{2},\text{ }t\in \left( \frac{1}{2},1\right]%
\end{array}%
\right. .
\end{equation*}

We require that the coefficients $-\frac{1}{2}\alpha \beta $, $\frac{1}{2}%
\left[ \left( \frac{1}{2}-\alpha \right) \left( \frac{1}{2}-\beta \right)
-\left( \frac{1}{2}-\gamma \right) \left( \frac{1}{2}-\delta \right) \right] 
$ and $\frac{1}{2}(1-\gamma )(1-\delta )$ be equal to zero. Hence, we
require that $\alpha =0$ or $\beta =0$ and $\gamma =1$ or $\delta =1$. If we
choose $\alpha =0$ and $\delta =1$ then we get $\beta +\gamma =1$. If we now
substitute $\alpha =0,$ $\gamma =1-\beta $ and $\delta =1$ in (\ref{j232a})
then we have

\begin{eqnarray}
&&\int\limits_{0}^{1}K_{2}(0,\beta ,1-\beta ,1,t)f^{\prime \prime }(t)dt
\label{j24a} \\
&=&-\int\limits_{0}^{1}K_{1}(0,\beta ,1-\beta ,1,t)f^{\prime }(t)dt  \notag
\\
&=&-\int\limits_{0}^{\frac{1}{2}}(t-\frac{\beta }{2})f^{\prime
}(t)dt-\int\limits_{\frac{1}{2}}^{1}(t-\frac{2-\beta }{2})f^{\prime }(t)dt 
\notag \\
&=&-\frac{\beta }{2}f(0)-(1-\beta )f(\frac{1}{2})-\frac{\beta }{2}%
f(1)+\int\limits_{0}^{1}f(t)dt.  \notag
\end{eqnarray}

We also have%
\begin{equation}
\left| \int\limits_{0}^{1}K_{2}(0,\beta ,1-\beta ,1,t)f^{\prime \prime
}(t)dt\right| \leq \left\| f^{\prime \prime }\right\| _{\infty
}\int\limits_{0}^{1}\left| K_{2}(0,\beta ,1-\beta ,1,t)\right| dt
\label{j24n}
\end{equation}%
and%
\begin{equation}
\int\limits_{0}^{1}\left| K_{2}(0,\beta ,1-\beta ,1,t)\right| dt=\frac{1}{2}%
\int\limits_{0}^{\frac{1}{2}}t\left| t-\beta \right| dt+\frac{1}{2}%
\int\limits_{\frac{1}{2}}^{1}\left| t-1+\beta \right| (1-t)dt.  \label{j25n}
\end{equation}%
We now define%
\begin{equation}
g(\beta )=\frac{1}{2}\int\limits_{0}^{\frac{1}{2}}t\left| t-\beta \right| dt+%
\frac{1}{2}\int\limits_{\frac{1}{2}}^{1}\left| t-1+\beta \right| (1-t)dt
\label{j26n}
\end{equation}%
and consider the problem%
\begin{equation}
\text{minimize }g(\beta )\text{, \ \ \ \ }\beta \in R.  \label{j27n}
\end{equation}%
Hence, we should like to find a global minimizer of $g$. Recall, a global
minimizer is a point $\beta ^{\ast }$ that satisfies%
\begin{equation}
g(\beta ^{\ast })\leq g(\beta )\text{, for all }\beta \in R.  \label{j28n}
\end{equation}

We consider the following cases:

(i) \ \ \ $\beta \leq 0$,

(ii) \ \ $0\leq \beta \leq \frac{1}{2}$,

(iii) \ $\beta \geq \frac{1}{2}$.

\textbf{The case (i).} If $\beta \leq 0$ then $t\left| t-\beta \right|
=t(t-\beta )$, for $t\in \left[ 0,\frac{1}{2}\right] $ and $\left| t-1+\beta
\right| \left| t-1\right| =(t-1+\beta )(t-1),$ for $t\in \left( \frac{1}{2},1%
\right] $. Thus,

\begin{eqnarray}
g(\beta ) &=&\frac{1}{2}\int\limits_{0}^{\frac{1}{2}}t(t-\beta )dt+\frac{1}{2%
}\int\limits_{\frac{1}{2}}^{1}(t-1+\beta )(t-1)dt  \label{jj10} \\
&=&\frac{1}{24}-\frac{\beta }{8}\geq \frac{1}{24}.  \notag
\end{eqnarray}

\textbf{The case (iii).} If $\beta \geq \frac{1}{2}$ then $t\left| t-\beta
\right| =t(\beta -t)$, for $t\in \left[ 0,\frac{1}{2}\right] $ and $\left|
t-1+\beta \right| \left| t-1\right| =(t-1+\beta )(1-t),$ for $t\in \left( 
\frac{1}{2},1\right] $. Thus,

\begin{eqnarray}
g(\beta ) &=&\frac{1}{2}\int\limits_{0}^{\frac{1}{2}}t(\beta -t)dt+\frac{1}{2%
}\int\limits_{\frac{1}{2}}^{1}(t-1+\beta )(1-t)dt  \label{jj11} \\
&=&\frac{\beta }{8}-\frac{1}{24}\geq \frac{1}{48}.  \notag
\end{eqnarray}

\textbf{The case (ii).} If $0\leq \beta \leq \frac{1}{2}$ then%
\begin{equation*}
t\left| t-\beta \right| =\left\{ 
\begin{array}{c}
t(\beta -t)\text{, \ \ \ \ }t\in \left[ 0,\beta \right] \\ 
t(t-\beta )\text{, \ \ }t\in \left( \beta ,\frac{1}{2}\right]%
\end{array}%
\right.
\end{equation*}%
and%
\begin{equation*}
\left| t-1+\beta \right| \left| t-1\right| =\left\{ 
\begin{array}{c}
(t-1+\beta )(t-1)\text{, \ \ }t\in \left[ \frac{1}{2},1-\beta \right] \\ 
(t-1+\beta )(1-t)\text{, \ \ \ \ }t\in \left( 1-\beta ,1\right]%
\end{array}%
\right. .
\end{equation*}%
Thus,%
\begin{eqnarray}
g(\beta ) &=&\frac{1}{2}\int\limits_{0}^{\beta }t(\beta -t)dt+\frac{1}{2}%
\int\limits_{\beta }^{\frac{1}{2}}t(t-\beta )dt  \label{jj12} \\
&&+\frac{1}{2}\int\limits_{\frac{1}{2}}^{1-\beta }(t-1+\beta )(t-1)dt+\frac{1%
}{2}\int\limits_{1-\beta }^{1}(t-1+\beta )(1-t)dt  \notag \\
&=&\frac{\beta ^{3}}{3}-\frac{\beta }{8}+\frac{1}{24}.  \notag
\end{eqnarray}%
We have%
\begin{equation}
g^{\prime }(\beta )=\beta ^{2}-\frac{1}{8}\text{ and }g^{\prime \prime
}(\beta )=2\beta .  \label{jj13}
\end{equation}%
From the equation $g^{\prime }(\beta )=0$ we find that $\beta _{1,2}=\pm 
\frac{\sqrt{2}}{4}$. Since $g^{\prime \prime }(\frac{\sqrt{2}}{4})>0$ we
conclude that $\beta =\frac{\sqrt{2}}{4}$ is, at least, a local minimizer.
We have

\begin{equation}
g(\frac{\sqrt{2}}{4})=\frac{2-\sqrt{2}}{48}.  \label{jj14}
\end{equation}%
From (\ref{jj10}), (\ref{jj11}) and (\ref{jj14}) we conclude that $\beta =%
\frac{\sqrt{2}}{4}$ is the global minimizer. If we now substitute $\beta =%
\frac{\sqrt{2}}{4}$ in (\ref{j24a}) then we get

\begin{eqnarray}
&&\int\limits_{0}^{1}K_{2}(0,\frac{\sqrt{2}}{4},1-\frac{\sqrt{2}}{4}%
,1,t)f^{\prime \prime }(t)dt  \label{jj15} \\
&=&\int\limits_{0}^{1}f(t)dt-\frac{\sqrt{2}}{8}f(0)-\left( 1-\frac{\sqrt{2}}{%
4}\right) f(\frac{1}{2})-\frac{\sqrt{2}}{8}f(1).  \notag
\end{eqnarray}%
The above quadrature formula is optimal in the sense described in Section 1.

From the previous considerations we can formulate the following result.

\begin{theorem}
\label{T1} Let $I\subset R$ be an open interval such that $\left[ 0,1\right]
\subset I$ \ and let $f:I\rightarrow R$ be a twice differentiable function
such that $f^{\prime \prime }$ is bounded and integrable. Then we have

\begin{equation}
\left| \int\limits_{0}^{1}f(t)dt-\frac{\sqrt{2}}{8}f(0)-\left( 1-\frac{\sqrt{%
2}}{4}\right) f(\frac{1}{2})-\frac{\sqrt{2}}{8}f(1)\right| \leq \frac{2-%
\sqrt{2}}{48}\left\| f^{\prime \prime }\right\| _{\infty }.  \label{jj16}
\end{equation}
\end{theorem}

\begin{remark}
\label{Rs} If we set $\beta =\frac{1}{3}$ in (\ref{j24a}) then we get the
well-known Simpson's rule:

\begin{equation}
\int\limits_{0}^{1}f(t)dt-\frac{1}{6}\left[ f(0)+4f(\frac{1}{2})+f(1)\right]
=\int\limits_{0}^{1}K_{2}(0,\frac{1}{3},\frac{2}{3},1,t)f^{\prime \prime
}(t)dt.  \label{jj17}
\end{equation}%
We have%
\begin{equation}
\left| \int\limits_{0}^{1}f(t)dt-\frac{1}{6}f(0)-\frac{2}{3}f(\frac{1}{2})-%
\frac{1}{6}f(1)\right| \leq \frac{\left\| f^{\prime \prime }\right\|
_{\infty }}{81}.  \label{jj18}
\end{equation}%
It is obvious that (\ref{jj16}) is a better estimate than (\ref{jj18}). Note
that (\ref{jj15}) and (\ref{jj17}) are 3-point quadrature rules of the same
(closed) type.
\end{remark}

If we consider the above problem on the interval $\left[ a,b\right] $ then
we get the following result.

\begin{theorem}
\label{T2ab} Let $I\subset R$ be an open interval such that $\left[ a,b%
\right] \subset I$ and let $f:I\rightarrow R$ be a twice differentiable
function such that $f^{\prime \prime }$ is bounded and integrable. Then we
have%
\begin{eqnarray}
&&\left| \int\limits_{a}^{b}f(t)dt-\left[ \frac{\sqrt{2}}{8}f(a)+\left( 1-%
\frac{\sqrt{2}}{4}\right) f(\frac{a+b}{2})+\frac{\sqrt{2}}{8}f(b)\right]
(b-a)\right|   \label{t2ab} \\
&\leq &\frac{2-\sqrt{2}}{48}\left\| f^{\prime \prime }\right\| _{\infty
}(b-a)^{3},  \notag
\end{eqnarray}%
where $\left\| f^{\prime \prime }\right\| _{\infty }=\underset{t\in \left[
a,b\right] }{\sup }\left| f^{\prime \prime }(t)\right| $.
\end{theorem}

\section{Error inequalities}

First we consider some basic properties of the spaces $L_{p}(a,b)$, for $%
p=1,2,\infty $. As we know, $X=(L_{2}(a,b),(\cdot ,\cdot ))$ is a Hilbert
space with the inner product%
\begin{equation}
(f,g)=\int_{a}^{b}f(t)g(t)dt.  \label{j1}
\end{equation}%
In the space $X$ the norm $\left\| \cdot \right\| _{2}$ is defined in the
usual way,%
\begin{equation}
\left\| f\right\| _{2}=\left( \int_{a}^{b}f(t)^{2}dt\right) ^{1/2}.
\label{j2}
\end{equation}%
We also consider the space $Y=(L_{2}(a,b),\left\langle \cdot ,\cdot
\right\rangle )$ where the inner product $\left\langle \cdot ,\cdot
\right\rangle $ is defined by%
\begin{equation}
\left\langle f,g\right\rangle =\frac{1}{b-a}\int_{a}^{b}f(t)g(t)dt.
\label{j3}
\end{equation}%
It is not difficult to see that $Y$ is a Hilbert space, too. In the space $Y$
the norm $\left\| \cdot \right\| $ is defined by%
\begin{equation}
\left\| f\right\| =\sqrt{\left\langle f,f\right\rangle }.  \label{j4}
\end{equation}%
We also define the Chebyshev functional%
\begin{equation}
T(f,g)=\left\langle f,g\right\rangle -\left\langle f,e\right\rangle
\left\langle g,e\right\rangle ,  \label{j5}
\end{equation}%
where $f,g\in L_{2}(a,b)$ and $e=1$. This functional satisfies the pre-Gr%
\"{u}ss inequality (\cite[p. 296]{MPF1}),%
\begin{equation}
T(f,g)^{2}\leq T(f,f)T(g,g).  \label{j6}
\end{equation}%
Specially, we define%
\begin{equation}
\sigma (f)=\sigma (f;a,b)=\sqrt{(b-a)T(f,f)}.  \label{j7}
\end{equation}%
The space $L_{1}(a,b)$ is a Banach space with the norm%
\begin{equation}
\left\| f\right\| _{1}=\int_{a}^{b}\left| f(t)\right| dt  \label{j8}
\end{equation}%
and the space $L_{\infty }(a,b)$ is also a Banach space with the norm%
\begin{equation}
\left\| f\right\| _{\infty }=\underset{t\in \left[ a,b\right] }{ess\sup }%
\left| f(t)\right| .  \label{j9}
\end{equation}%
If $f\in L_{1}(a,b)$ and $g\in L_{\infty }(a,b)$ then we have%
\begin{equation}
\left| (f,g)\right| \leq \left\| f\right\| _{1}\left\| g\right\| _{\infty }.
\label{j10}
\end{equation}%
More about the above mentioned spaces can be found, for example, in \cite%
{AT1}.

Finally, we define the functional%
\begin{eqnarray}
Q(f) &=&Q(f;a,b)  \label{j11} \\
&=&\int_{a}^{b}f(t)dt-\left[ \frac{\sqrt{2}}{8}f(a)+\left( 1-\frac{\sqrt{2}}{%
4}\right) f(\frac{a+b}{2})+\frac{\sqrt{2}}{8}f(b)\right] (b-a).  \notag
\end{eqnarray}%
We also need the following lemma.

\begin{lemma}
\label{L1} Let%
\begin{equation}
f(t)=\left\{ 
\begin{array}{c}
f_{1}(t),\text{ \ }t\in \left[ a,x_{0}\right] \\ 
f_{2}(t),\text{ \ }t\in \left( x_{0},b\right]%
\end{array}%
\right. ,  \label{a1}
\end{equation}%
where $x_{0}\in \left[ a,b\right] $, $f_{1}\in C^{1}(a,x_{0})$, $f_{2}\in
C^{1}(x_{0},b)$. If $f_{1}(x_{0})=f_{2}(x_{0})$ then $f$ is an absolutely
continuous function.
\end{lemma}

A proof of this lemma can be found in \cite{U2}.

\begin{theorem}
Let $f:\left[ 0,1\right] \rightarrow R$ be an absolutely continuous function
such that $f^{\prime }\in L_{1}(0,1)$ and there exist real numbers $\gamma
_{1},\Gamma _{1}$ such that $\gamma _{1}\leq f^{\prime }(t)\leq \Gamma _{1}$%
, $t\in \left[ 0,1\right] $. Then%
\begin{equation}
\left| Q(f;0,1)\right| \leq \frac{\Gamma _{1}-\gamma _{1}}{32}(5-2\sqrt{2}),
\label{j14}
\end{equation}%
\begin{equation}
\left| Q(f;0,1)\right| \leq \left( \frac{1}{2}-\frac{\sqrt{2}}{8}\right)
(S-\gamma _{1}),  \label{j15}
\end{equation}%
\begin{equation}
\left| Q(f;0,1)\right| \leq \left( \frac{1}{2}-\frac{\sqrt{2}}{8}\right)
(\Gamma _{1}-S),  \label{j16}
\end{equation}%
where $Q(f;0,1)$ is defined by (\ref{j11}) and $S=f(1)-f(0)$.
\end{theorem}

\begin{proof}
We define the function%
\begin{equation}
p_{1}(t)=\left\{ 
\begin{array}{c}
t-\frac{\sqrt{2}}{8}\text{, \ \ \ \ \ }t\in \left[ 0,\frac{1}{2}\right] \\ 
t-1+\frac{\sqrt{2}}{8}\text{, }t\in \left( \frac{1}{2},1\right]%
\end{array}%
\right. .  \label{j17}
\end{equation}%
It is easy to verify that%
\begin{equation}
(p_{1},f^{\prime })=-Q(f;0,1).  \label{j18}
\end{equation}%
On the other hand, we have%
\begin{equation}
\left( f^{\prime }-\frac{\Gamma _{1}+\gamma _{1}}{2},p_{1}\right)
=(f^{\prime },p_{1}),  \label{j19}
\end{equation}%
since $(p_{1},e)=0$. From (\ref{j10}) we get%
\begin{equation}
\left| \left( f^{\prime }-\frac{\Gamma _{1}+\gamma _{1}}{2},p_{1}\right)
\right| \leq \left\| f^{\prime }-\frac{\Gamma _{1}+\gamma _{1}}{2}\right\|
_{\infty }\left\| p_{1}\right\| _{1}\leq \frac{\Gamma _{1}-\gamma _{1}}{32}%
(5-2\sqrt{2}),  \label{j20}
\end{equation}%
since%
\begin{equation*}
\left\| f^{\prime }-\frac{\Gamma _{1}+\gamma _{1}}{2}\right\| _{\infty }\leq 
\frac{\Gamma _{1}-\gamma _{1}}{2}
\end{equation*}%
and%
\begin{equation*}
\left\| p_{1}\right\| _{1}=\frac{5}{16}-\frac{\sqrt{2}}{8}.
\end{equation*}%
From (\ref{j18})-(\ref{j20}) we see that (\ref{j14}) holds. We now prove
that (\ref{j15}) holds. We have%
\begin{equation*}
\left| (f^{\prime }-\gamma _{1},p_{1})\right| \leq \left\| p_{1}\right\|
_{\infty }\left\| f^{\prime }-\gamma _{1}\right\| _{1}=\left( \frac{1}{2}-%
\frac{\sqrt{2}}{8}\right) (S-\gamma _{1}),
\end{equation*}%
since%
\begin{equation*}
\left\| p_{1}\right\| _{\infty }=\frac{1}{2}-\frac{\sqrt{2}}{8}
\end{equation*}%
and%
\begin{equation*}
\left\| f^{\prime }-\gamma _{1}\right\| _{1}=\int_{0}^{1}(f^{\prime
}(t)-\gamma _{1})dt=f(1)-f(0)-\gamma _{1}.
\end{equation*}%
In a similar way we can prove that (\ref{j16}) holds.
\end{proof}

\begin{remark}
\label{R1} Note that we can apply the estimate (\ref{j14}) only if the first
derivative $f^{\prime }$ is bounded. It means that we cannot use (\ref{j14})
to estimate directly the error when approximating the integral of such a
well-behaved function as $f(t)=\sqrt{t}$ on $\left[ 0,1\right] $, (since $%
f^{\prime }(t)=1/(2\sqrt{t})$ is unbounded on $\left[ 0,1\right] $). On the
other hand, we can use the estimation (\ref{j15}), (since $\gamma =1/2$ on $%
\left[ 0,1\right] $ for the given function).
\end{remark}

\begin{theorem}
\label{T2} Let $f:\left[ a,b\right] \rightarrow R$ be an absolutely
continuous function such that $f^{\prime }\in L_{1}(a,b)$ and there exist
real numbers $\gamma _{1},\Gamma _{1}$ such that $\gamma _{1}\leq f^{\prime
}(t)\leq \Gamma _{1}$, $t\in \left[ a,b\right] $. Then%
\begin{equation}
\left| Q(f;a,b)\right| \leq \frac{\Gamma _{1}-\gamma _{1}}{32}(5-2\sqrt{2}%
)(b-a)^{2},  \label{j21}
\end{equation}%
\begin{equation}
\left| Q(f;a,b)\right| \leq \left( \frac{1}{2}-\frac{\sqrt{2}}{8}\right)
(S-\gamma _{1})(b-a)^{2},  \label{j22}
\end{equation}%
\begin{equation}
\left| Q(f;a,b)\right| \leq \left( \frac{1}{2}-\frac{\sqrt{2}}{8}\right)
(\Gamma _{1}-S)(b-a)^{2},  \label{j23}
\end{equation}%
where $Q(f;a,b)$ is defined by (\ref{j11}) and $S=(f(b)-f(a))/(b-a)$.
\end{theorem}

\begin{theorem}
\label{T3} Let $f:\left[ 0,1\right] \rightarrow R$ be an absolutely
continuous function such that $f^{\prime }\in L_{2}(0,1)$. Then%
\begin{equation}
\left| Q(f;0,1)\right| \leq \sqrt{\frac{11}{96}-\frac{\sqrt{2}}{16}}\sigma
(f^{\prime };0,1),  \label{j24}
\end{equation}%
where $\sigma (f;0,1)$ is defined by (\ref{j7}). The inequality (\ref{j24})
is sharp in the sense that the constant $\sqrt{\frac{11}{96}-\frac{\sqrt{2}}{%
16}}$ cannot be replaced by a smaller one.
\end{theorem}

\begin{proof}
Let $p_{1}$ be defined by (\ref{j17}). We have%
\begin{equation*}
\left\langle p_{1},f^{\prime }\right\rangle =-Q(f;0,1),
\end{equation*}%
since (\ref{j18}) holds and $\left\langle f,g\right\rangle =(f,g)$ if $\left[
a,b\right] =\left[ 0,1\right] $. On the other hand, we have%
\begin{equation*}
\left\langle p_{1},f^{\prime }\right\rangle =T(f^{\prime },p_{1}),
\end{equation*}%
since $\left\langle p_{1},e\right\rangle =0$. From (\ref{j6}) it follows%
\begin{eqnarray*}
\left| T(f^{\prime },p_{1})\right| &\leq &\sqrt{T(p_{1},p_{1})}\sqrt{%
T(f^{\prime },f^{\prime })}=\left\| p_{1}\right\| \sigma (f^{\prime };0,1) \\
&=&\sqrt{\frac{11}{96}-\frac{\sqrt{2}}{16}}\sigma (f^{\prime };0,1),
\end{eqnarray*}%
since%
\begin{equation*}
\left\| p_{1}\right\| =\sqrt{\frac{11}{96}-\frac{\sqrt{2}}{16}}.
\end{equation*}%
Hence, the inequality (\ref{j24}) is proved. We have to prove that this
inequality is sharp. For that purpose, we define the function%
\begin{equation}
f(t)=\left\{ 
\begin{array}{c}
\frac{1}{2}t^{2}-\frac{\sqrt{2}}{8}t\text{, \ \ }t\in \left[ 0,\frac{1}{2}%
\right] \\ 
\frac{1}{2}t^{2}-(1-\frac{\sqrt{2}}{8})t+\frac{1}{2}-\frac{\sqrt{2}}{8}\text{%
, }t\in \left( \frac{1}{2},1\right]%
\end{array}%
\right.  \label{j25}
\end{equation}%
such that $f^{\prime }(t)=p_{1}(t)$. From Lemma \ref{L1} we see that the
function $f$, defined by (\ref{j25}), is an absolutely continuous function.
For this function the left-hand side of (\ref{j24}) becomes%
\begin{equation*}
L.H.S.(\ref{j24})=\left| -\frac{11}{96}+\frac{\sqrt{2}}{16}\right| .
\end{equation*}%
The right-hand side of (\ref{j24}) becomes%
\begin{equation*}
R.H.S.(\ref{j24})=\frac{11}{96}-\frac{\sqrt{2}}{16}.
\end{equation*}%
We see that $L.H.S.(\ref{j24})=R.H.S.(\ref{j24})$. Thus, (\ref{j24}) is
sharp.
\end{proof}

\begin{remark}
\label{R3} The estimate (\ref{j14}) is better than the estimate (\ref{j24}).
However, note that the estimate (\ref{j14}) can be applied only if $%
f^{\prime }$ is bounded. On the other hand, the estimate (\ref{j1}) can be
applied for an absolutely continuous function if $f^{\prime }\in L_{2}(a,b)$.
\end{remark}

There are many examples where we cannot apply the estimate (\ref{j14}) but
we can apply (\ref{j24}).

\begin{example}
\label{E1}Let us consider the integral $\int\limits_{0}^{1}\sqrt[3]{\sin
t^{2}}dt.$ We have%
\begin{equation*}
f(t)=\sqrt[3]{\sin t^{2}}\text{ and }f^{\prime }(t)=\frac{2t\cos t^{2}}{3%
\sqrt[3]{\sin ^{2}t^{2}}}
\end{equation*}%
such that $f^{\prime }(t)\rightarrow \infty $, \ $t\rightarrow 0$ and we
cannot apply the estimate (\ref{j14}). On the other hand, we have%
\begin{equation*}
\int\limits_{0}^{1}\left[ f^{\prime }(t)\right] ^{2}dt\leq \frac{4}{9}%
\underset{t\in \left[ 0,1\right] }{\max }\frac{t^{2}\cos t^{2}}{\sin t^{2}}%
\int\limits_{0}^{1}\frac{dt}{\sqrt[3]{\sin t^{2}}}\leq \frac{16}{9},
\end{equation*}%
i.e. $\left\| f^{\prime }\right\| _{2}\leq \frac{4}{3}$ and we can apply the
estimate (\ref{j24}).
\end{example}

\begin{theorem}
\label{T4} Let $f:\left[ a,b\right] \rightarrow R$ be an absolutely
continuous function such that $f^{\prime }\in L_{2}(a,b)$. Then%
\begin{equation}
\left| Q(f;a,b)\right| \leq \sqrt{\frac{11}{96}-\frac{\sqrt{2}}{16}}\sigma
(f^{\prime };a,b)(b-a)^{3/2},  \label{j28}
\end{equation}%
where $\sigma (f;a,b)$ is defined by (\ref{j7}). The inequality (\ref{j28})
is sharp in the sense that the constant $\sqrt{\frac{11}{96}-\frac{\sqrt{2}}{%
16}}$ cannot be replaced by a smaller one.
\end{theorem}

We define 
\begin{equation}
P(f;a,b)=\frac{(b-a)^{2}}{96}\left( 4-3\sqrt{2}\right) \left[ f^{\prime
}(b)-f^{\prime }(a)\right] .  \label{k0}
\end{equation}

\begin{theorem}
\label{T5} Let $f^{\prime }:\left[ 0,1\right] \rightarrow R$ be an
absolutely continuous function such that $f^{\prime \prime }\in L_{1}(0,1)$
and there exist real numbers $\gamma _{2},\Gamma _{2}$ such that $\gamma
_{2}\leq f^{\prime \prime }(t)\leq \Gamma _{2}$, $t\in \left[ 0,1\right] $.
Then%
\begin{equation}
\left| Q(f;0,1)-P(f;0,1)\right| \leq \frac{\Gamma _{2}-\gamma _{2}}{2}\left( 
\frac{5}{96}\sqrt{6}-\frac{29}{432}\sqrt{3}\right) ,  \label{k1}
\end{equation}%
\begin{equation}
\left| Q(f;0,1)-P(f;0,1)\right| \leq \left( \frac{1}{12}-\frac{\sqrt{2}}{32}%
\right) (S_{1}-\gamma _{2}),  \label{k2}
\end{equation}%
\begin{equation}
\left| Q(f;0,1)-P(f;0,1)\right| \leq \left( \frac{1}{12}-\frac{\sqrt{2}}{32}%
\right) (\Gamma _{2}-S_{1}),  \label{k3}
\end{equation}%
where $Q(f;0,1)$ and $P(f;0,1)$ are defined by (\ref{j11}) and (\ref{k0}),
respectively and $S_{1}=f^{\prime }(1)-f^{\prime }(0)$.
\end{theorem}

\begin{proof}
We define the function%
\begin{equation}
\tilde{p}_{2}(t)=\left\{ 
\begin{array}{c}
\frac{1}{2}t(t-\frac{\sqrt{2}}{4})+\frac{\sqrt{2}}{32}-\frac{1}{24}\text{, \
\ \ \ \ \ \ \ \ \ \ \ \ \ \ \ }t\in \left[ 0,\frac{1}{2}\right] \\ 
\frac{1}{2}(t-1)(t-1+\frac{\sqrt{2}}{4})+\frac{\sqrt{2}}{32}-\frac{1}{24}%
\text{, }t\in \left( \frac{1}{2},1\right]%
\end{array}%
\right. .  \label{k4}
\end{equation}%
Let $p_{1}$ be defined by (\ref{j17}). Then we have%
\begin{equation}
(\tilde{p}_{2},f^{\prime \prime })=-(p_{1},f^{\prime
})-P(f;0,1)=Q(f;0,1)-P(f;0,1)  \label{k5}
\end{equation}%
since (\ref{j18}) holds.

On the other hand, we have%
\begin{equation}
\left( f^{\prime \prime }-\frac{\Gamma _{2}+\gamma _{2}}{2},\tilde{p}%
_{2}\right) =(f^{\prime \prime },\tilde{p}_{2}),  \label{k6}
\end{equation}%
since $(\tilde{p}_{2},e)=0$. From (\ref{j6}) we get%
\begin{eqnarray}
\left| \left( f^{\prime \prime }-\frac{\Gamma _{2}+\gamma _{2}}{2},\tilde{p}%
_{2}\right) \right| &\leq &\left\| f^{\prime \prime }-\frac{\Gamma
_{2}+\gamma _{2}}{2}\right\| _{\infty }\left\| \tilde{p}_{2}\right\| _{1}
\label{k7} \\
&\leq &\left( \frac{5}{96}\sqrt{6}-\frac{29}{432}\sqrt{3}\right) \frac{%
\Gamma _{2}-\gamma _{2}}{2},  \notag
\end{eqnarray}%
since%
\begin{equation*}
\left\| f^{\prime \prime }-\frac{\Gamma _{2}+\gamma _{2}}{2}\right\|
_{\infty }\leq \frac{\Gamma _{2}-\gamma _{2}}{2}
\end{equation*}%
and%
\begin{equation*}
\left\| \tilde{p}_{2}\right\| _{1}=\frac{5}{96}\sqrt{6}-\frac{29}{432}\sqrt{3%
}.
\end{equation*}%
From (\ref{k5})-(\ref{k7}) we see that (\ref{k1}) holds.

We now prove that (\ref{k2}) holds. We have%
\begin{equation*}
\left| \left( f^{\prime \prime }-\gamma _{2},\tilde{p}_{2}\right) \right|
\leq \left\| f^{\prime \prime }-\gamma _{2}\right\| _{1}\left\| \tilde{p}%
_{2}\right\| _{\infty }=\left( \frac{1}{12}-\frac{\sqrt{2}}{32}\right)
(S_{1}-\gamma _{2}),
\end{equation*}%
since%
\begin{equation*}
\left\| f^{\prime \prime }-\gamma _{2}\right\| _{1}=\int_{0}^{1}(f^{\prime
\prime }(t)-\gamma _{2})dt=f^{\prime }(1)-f^{\prime }(0)-\gamma _{2}
\end{equation*}%
and%
\begin{equation*}
\left\| \tilde{p}_{2}\right\| _{\infty }=\frac{1}{12}-\frac{\sqrt{2}}{32}.
\end{equation*}%
In a similar way we can prove that (\ref{k3}) holds.
\end{proof}

\begin{remark}
\label{R4} Note that we can apply the estimate (\ref{k1}) only if the second
derivative $f^{\prime \prime }$ is bounded. It means that we cannot use (\ref%
{k1}) to estimate directly the error when approximating the integral of such
a well-behaved function as $f(t)=\sqrt{t^{3}}$ on $\left[ 0,1\right] $,
(since $f^{\prime \prime }(t)=3/(4\sqrt{t})$ is unbounded on $\left[ 0,1%
\right] $). On the other hand, we can use the estimation (\ref{k2}), (since $%
\gamma =3/4$ on $\left[ 0,1\right] $ for the given function).
\end{remark}

\begin{theorem}
\label{T6} Let $f^{\prime }:\left[ a,b\right] \rightarrow R$ be an
absolutely continuous function such that $f^{\prime \prime }\in L_{1}(a,b)$
and there exist real numbers $\gamma _{2},\Gamma _{2}$ such that $\gamma
_{2}\leq f^{\prime \prime }(t)\leq \Gamma _{2}$, $t\in \left[ a,b\right] $.
Then%
\begin{equation}
\left| Q(f;a,b)-P(f;a,b)\right| \leq \frac{\Gamma _{2}-\gamma _{2}}{2}\left( 
\frac{5}{96}\sqrt{6}-\frac{29}{432}\sqrt{3}\right) (b-a)^{3},  \label{k11}
\end{equation}%
\begin{equation}
\left| Q(f;a,b)-P(f;a,b)\right| \leq \left( \frac{1}{12}-\frac{\sqrt{2}}{32}%
\right) (S_{1}-\gamma _{2})(b-a)^{3},  \label{k12}
\end{equation}%
\begin{equation}
\left| Q(f;a,b)-P(f;a,b)\right| \leq \left( \frac{1}{12}-\frac{\sqrt{2}}{32}%
\right) (\Gamma _{2}-S_{1})(b-a)^{3},  \label{k13}
\end{equation}%
where $Q(f;a,b)$ and $P(f;a,b)$ are defined by (\ref{j11}) and (\ref{k0}),
respectively and $S_{1}=(f^{\prime }(b)-f^{\prime }(a))/(b-a)$.
\end{theorem}

\begin{theorem}
\label{T7} Let $f^{\prime }:\left[ 0,1\right] \rightarrow R$ be an
absolutely continuous function such that $f^{\prime \prime }\in L_{2}(0,1)$.
Then%
\begin{equation}
\left| Q(f;0,1)-P(f;0,1)\right| \leq \sqrt{\frac{47}{23040}-\frac{\sqrt{2}}{%
768}}\sigma (f^{\prime \prime };0,1)  \label{k10}
\end{equation}%
where $\sigma (f;0,1)$ is defined by (\ref{j7}). The inequality (\ref{k10})
is sharp in the sense that the constant $\sqrt{\frac{47}{23040}-\frac{\sqrt{2%
}}{768}}$ cannot be replaced by a smaller one.
\end{theorem}

\begin{proof}
We define the function%
\begin{equation}
p_{2}(t)=\left\{ 
\begin{array}{c}
\frac{1}{2}t(t-\frac{\sqrt{2}}{4})\text{, \ \ \ \ \ \ \ \ \ \ \ \ \ \ \ }%
t\in \left[ 0,\frac{1}{2}\right] \\ 
\frac{1}{2}(t-1)(t-1+\frac{\sqrt{2}}{4})\text{, }t\in \left( \frac{1}{2},1%
\right]%
\end{array}%
\right. .  \label{k14}
\end{equation}%
Then we have%
\begin{equation}
\left\langle \tilde{p}_{2},f^{\prime \prime }\right\rangle =\left\langle
p_{2},f^{\prime \prime }\right\rangle -\left\langle p_{2},e\right\rangle
\left\langle f^{\prime \prime },e\right\rangle  \label{k15}
\end{equation}%
since $\tilde{p}_{2}=p_{2}-\left\langle p_{2},e\right\rangle $. From (\ref%
{k5}) and (\ref{k15}) it follows%
\begin{equation}
T(p_{2},f^{\prime \prime })=Q(f;0,1)-P(f;0,1),  \label{k16}
\end{equation}%
since $\left\langle \tilde{p}_{2},f^{\prime \prime }\right\rangle =(\tilde{p}%
_{2},f^{\prime \prime })$ if $\left[ a,b\right] =\left[ 0,1\right] $. From (%
\ref{j6}) we get%
\begin{equation}
\left| T(p_{2},f^{\prime \prime })\right| \leq \sqrt{T(p_{2},p_{2})}\sqrt{%
T(f^{\prime \prime },f^{\prime \prime })}=\sqrt{\frac{47}{23040}-\frac{\sqrt{%
2}}{768}}\sigma (f^{\prime \prime };0,1),  \label{k17}
\end{equation}%
since%
\begin{equation*}
T(p_{2},p_{2})=\frac{47}{23040}-\frac{\sqrt{2}}{768}.
\end{equation*}%
From (\ref{k16}) and (\ref{k17}) we see that (\ref{k10}) holds.

We now prove that (\ref{k10}) is sharp. For that purpose we define the
function%
\begin{equation}
f(t)=\left\{ 
\begin{array}{c}
\frac{t^{4}}{24}-\frac{\sqrt{2}}{48}t^{3}\text{, \ \ \ \ \ \ \ \ \ \ \ \ \ \
\ \ \ \ \ \ \ \ \ \ \ \ \ \ \ \ \ \ \ \ \ \ \ \ \ \ \ \ \ \ \ \ \ \ \ \ \ \
\ \ \ \ \ \ \ \ \ \ \ \ \ \ \ \ \ \ \ \ }t\in \left[ 0,\frac{1}{2}\right] \\ 
\frac{t^{4}}{24}-(\frac{1}{6}-\frac{\sqrt{2}}{48})t^{3}+(\frac{1}{4}-\frac{%
\sqrt{2}}{16})t^{2}-(\frac{1}{8}-\frac{\sqrt{2}}{32})t+\frac{1}{48}-\frac{%
\sqrt{2}}{192}\text{, }t\in \left( \frac{1}{2},1\right]%
\end{array}%
\right.  \label{k19}
\end{equation}%
such that%
\begin{equation}
f^{\prime }(t)=\left\{ 
\begin{array}{c}
\frac{t^{3}}{6}-\frac{\sqrt{2}}{16}t^{2}\text{, \ \ \ \ \ \ \ \ \ \ \ \ \ \
\ \ \ \ \ \ \ \ \ \ \ \ \ \ \ \ \ \ \ \ \ \ \ \ \ \ \ \ \ \ \ \ \ \ \ \ \ }%
t\in \left[ 0,\frac{1}{2}\right] \\ 
\frac{t^{3}}{6}-(\frac{1}{2}-\frac{\sqrt{2}}{16})t^{2}+(\frac{1}{2}-\frac{%
\sqrt{2}}{8})t-(\frac{1}{8}-\frac{\sqrt{2}}{32})\text{, }t\in \left( \frac{1%
}{2},1\right]%
\end{array}%
\right.  \label{k19a}
\end{equation}

and $f^{\prime \prime }(t)=p_{2}(t)$. From Lemma \ref{L1} we see that the
function $f^{\prime }$, defined by (\ref{k19a}), is an absolutely continuous
function. For the function defined by (\ref{k19}) the left-hand side of (\ref%
{k10}) becomes%
\begin{equation*}
L.H.S.(\ref{k10})=\frac{47}{23040}-\frac{\sqrt{2}}{768}.
\end{equation*}%
The right-hand side of (\ref{k10}) becomes%
\begin{equation*}
R.H.S.(\ref{k10})=\frac{47}{23040}-\frac{\sqrt{2}}{768}.
\end{equation*}%
We see that $L.H.S.(\ref{k10})=R.H.S.(\ref{k10})$. Thus, (\ref{k10}) is
sharp.
\end{proof}

\begin{remark}
\label{R5} The estimation (\ref{k1}) is better than the estimation (\ref{k10}%
). However, note that we can apply the estimate (\ref{k1}) only if the
second derivative $f^{\prime \prime }$ is bounded. It means that we cannot
use (\ref{k1}) to estimate directly the error when approximating the
integral of such a well-behaved function as $f(t)=\sqrt[3]{t^{5}}$ on $\left[
0,1\right] $, (since $f^{\prime \prime }(t)=10/(9\sqrt[3]{t})$ is unbounded
on $\left[ 0,1\right] $). On the other hand, we can use the estimation (\ref%
{k10}), (since $\left\| f^{\prime \prime }\right\| _{2}^{2}=\frac{100}{27}$
for the given function).
\end{remark}

\begin{theorem}
\label{T8} Let $f^{\prime }:\left[ a,b\right] \rightarrow R$ be an
absolutely continuous function such that $f^{\prime \prime }\in L_{2}(a,b)$.
Then%
\begin{equation}
\left| Q(f;a,b)-P(f;a,b)\right| \leq \sqrt{\frac{47}{23040}-\frac{\sqrt{2}}{%
768}}\sigma (f^{\prime \prime };a,b)(b-a)^{5/2},  \label{k20}
\end{equation}%
where $\sigma (f;a,b)$ is defined by (\ref{j7}). The inequality (\ref{k20})
is sharp in the sense that the constant $\sqrt{\frac{47}{23040}-\frac{\sqrt{2%
}}{768}}$ cannot be replaced by a smaller one.
\end{theorem}

\section{Applications in numerical integration}

Let $\pi =\left\{ x_{0}=a<x_{1}<\cdots <x_{n}=b\right\} $ be a given
subdivision of the interval $\left[ a,b\right] $ such that $%
h_{i}=x_{i+1}-x_{i}=h=(b-a)/n$. From (\ref{j11}) we get

\begin{eqnarray*}
&&Q(f;x_{i},x_{i+1}) \\
&=&\int_{x_{i}}^{x_{i+1}}f(t)dt-\left[ \frac{\sqrt{2}}{8}f(x_{i})+\left( 1-%
\frac{\sqrt{2}}{4}\right) f(\frac{x_{i}+x_{i+1}}{2})+\frac{\sqrt{2}}{8}%
f(x_{i+1})\right] h.
\end{eqnarray*}

If we now sum the above relation over $i$ from $0$ to $n-1$ then we get

\begin{eqnarray*}
&&\sum\limits_{i=0}^{n-1}Q(f;x_{i},x_{i+1}) \\
&=&\int_{a}^{b}f(t)dt-\frac{\sqrt{2}h}{8}\left[ f(a)+f(b)\right] \\
&&-\frac{\sqrt{2}h}{4}\sum\limits_{i=1}^{n-1}f(x_{i})-\left( 1-\frac{\sqrt{2}%
}{4}\right) h\sum\limits_{i=1}^{n-1}f(\frac{x_{i}+x_{i+1}}{2}).
\end{eqnarray*}

We introduce the notation

\begin{equation}
S(f;a,b)=\sum\limits_{i=0}^{n-1}Q(f;x_{i},x_{i+1}).  \label{g0}
\end{equation}%
We also define%
\begin{equation}
P_{n}(f;a,b)=\frac{(b-a)^{2}}{96n^{2}}\left( 4-3\sqrt{2}\right) \left[
f^{\prime }(b)-f^{\prime }(a)\right] ,  \label{k01}
\end{equation}%
\begin{equation}
\sigma _{n}(f)=\sum\limits_{i=0}^{n-1}\sqrt{\frac{b-a}{n}\left\| f^{\prime
}\right\| _{2}^{2}-\left[ f(x_{i+1})-f(x_{i})\right] ^{2}}  \label{g1}
\end{equation}%
and%
\begin{equation}
\omega _{n}(f)=\left[ (b-a)\left\| f^{\prime }\right\| _{2}^{2}-\frac{1}{n}%
(f(b)-f(a))^{2}\right] ^{1/2}.  \label{g1a}
\end{equation}

\begin{theorem}
\label{T4ab} Under the assumptions of Theorem \ref{T2ab} we have%
\begin{eqnarray*}
&&\left| \int\limits_{a}^{b}f(t)dt-\sum\limits_{i=0}^{n-1}\left[ \frac{\sqrt{%
2}}{8}f(x_{i})+\left( 1-\frac{\sqrt{2}}{4}\right) f(\frac{x_{i}+x_{i+1}}{2})+%
\frac{\sqrt{2}}{8}f(x_{i+1})\right] h\right|  \\
&\leq &\frac{2-\sqrt{2}}{48n^{2}}\left\| f^{\prime \prime }\right\| _{\infty
}(b-a)^{3},
\end{eqnarray*}%
where $\left\{ a=x_{0}<x_{1}<\cdots <x_{n}=b\right\} $ is a uniform
subdivision of $\left[ a,b\right] $, i.e. $x_{i}=a+ih,$ $h=(b-a)/n$, $%
i=0,1,...,n$.
\end{theorem}

\begin{proof}
Apply Theorem \ref{T2ab} to the intervals $\left[ x_{i},x_{i+1}\right] $ and
sum.
\end{proof}

\begin{theorem}
\label{T1p} Under the assumptions of Theorem \ref{T2} we have%
\begin{equation*}
\left| S(f;a,b\right| \leq \frac{\Gamma _{1}-\gamma _{1}}{32n}(5-2\sqrt{2}%
)(b-a)^{2},
\end{equation*}%
\begin{equation*}
\left| S(f;a,b\right| \leq \frac{S-\gamma _{1}}{n}\left( \frac{1}{2}-\frac{%
\sqrt{2}}{8}\right) (b-a)^{2},
\end{equation*}%
\begin{equation*}
\left| S(f;a,b\right| \leq \frac{\Gamma _{1}-S}{n}\left( \frac{1}{2}-\frac{%
\sqrt{2}}{8}\right) (b-a)^{2},
\end{equation*}%
where $S(f;a,b)$ is defined by (\ref{g0}) and $\left\{ a=x_{0}<x_{1}<\cdots
<x_{n}=b\right\} $ is a uniform subdivision of $\left[ a,b\right] $, i.e. $%
x_{i}=a+ih,$ $h=(b-a)/n$, $i=0,1,...,n$.
\end{theorem}

\begin{proof}
Apply Theorem \ref{T2} to the intervals $\left[ x_{i},x_{i+1}\right] $ and
sum. Note that%
\begin{equation*}
\sum\limits_{i=0}^{n-1}\left[ f(x_{i+1})-f(x_{i})\right] =f(b)-f(a).
\end{equation*}
\end{proof}

\begin{theorem}
\label{T2p} Under the assumptions of Theorem \ref{T4} we have%
\begin{equation}
\left| S(f;a,b\right| \leq \sqrt{\frac{11}{96}-\frac{\sqrt{2}}{16}}\frac{b-a%
}{n}\sigma _{n}(f)\leq \sqrt{\frac{11}{96}-\frac{\sqrt{2}}{16}}\frac{b-a}{%
\sqrt{n}}\omega _{n}(f),  \label{g2}
\end{equation}%
where $S(f;a,b)$, $\sigma _{n}(f)$ and $\omega _{n}(f)$ are defined by (\ref%
{g0}), (\ref{g1}) and (\ref{g1a}), respectively and $\left\{
a=x_{0}<x_{1}<\cdots <x_{n}=b\right\} $ is a uniform subdivision of $\left[
a,b\right] $, i.e. $x_{i}=a+ih,$ $h=(b-a)/n$, $i=0,1,...,n$.
\end{theorem}

\begin{proof}
We apply Theorem \ref{T4} to the interval $\left[ x_{i},x_{i+1}\right] $ and
sum. Then we have%
\begin{eqnarray*}
&&\left| S(f;a,b)\right| \\
&\leq &\sqrt{\frac{11}{96}-\frac{\sqrt{2}}{16}}h^{3/2}\sum\limits_{i=0}^{n-1}%
\left[ \left\| f^{\prime }\right\| _{2}^{2}-\frac{1}{h}\left(
f(x_{i+1})-f(x_{i})\right) ^{2}\right] ^{1/2}.
\end{eqnarray*}%
From the above relation and the fact $h=(b-a)/n$ we see that the first
inequality in (\ref{g2}) holds.

Using the Cauchy inequality we get%
\begin{eqnarray}
&&\sum\limits_{i=0}^{n-1}\left[ \left\| f^{\prime }\right\| _{2}^{2}-\frac{1%
}{h}\left( f(x_{i+1})-f(x_{i})\right) ^{2}\right] ^{1/2}  \label{g7a} \\
&\leq &n\left[ \left\| f^{\prime }\right\| _{2}^{2}-\frac{1}{b-a}%
\sum\limits_{i=0}^{n-1}\left( f(x_{i+1})-f(x_{i})\right) ^{2}\right] ^{1/2} 
\notag \\
&\leq &n\left[ \left\| f^{\prime }\right\| _{2}^{2}-\frac{1}{b-a}\frac{1}{n}%
\left( f(b)-f(a)\right) ^{2}\right] ^{1/2}.  \notag
\end{eqnarray}%
Thus the second inequality in (\ref{g2}) holds, too.
\end{proof}

\begin{theorem}
\label{T3p} Under the assumptions of Theorem \ref{T6} we have%
\begin{equation*}
\left| S(f;a,b)-P_{n}(f;a,b)\right| \leq \frac{\Gamma _{2}-\gamma _{2}}{2n}%
\left( \frac{5}{96}\sqrt{6}-\frac{29}{432}\sqrt{3}\right) (b-a)^{3},
\end{equation*}%
\begin{equation*}
\left| S(f;a,b)-P_{n}(f;a,b)\right| \leq \left( \frac{1}{12}-\frac{\sqrt{2}}{%
32}\right) \frac{(S_{1}-\gamma _{2})}{n}(b-a)^{3},
\end{equation*}%
\begin{equation*}
\left| S(f;a,b)-P_{n}(f;a,b)\right| \leq \left( \frac{1}{12}-\frac{\sqrt{2}}{%
32}\right) \frac{(\Gamma _{2}-S_{1})}{n}(b-a)^{3},
\end{equation*}%
where $S(f;a,b)$ and $P_{n}(f;a,b)$ are defined by (\ref{g0}) and (\ref{k01}%
), respectively and $\left\{ a=x_{0}<x_{1}<\cdots <x_{n}=b\right\} $ is a
uniform subdivision of $\left[ a,b\right] $, i.e. $x_{i}=a+ih,$ $h=(b-a)/n$, 
$i=0,1,...,n$.
\end{theorem}

\begin{proof}
The proof of this theorem is similar to the proof of Theorem \ref{T1p}. Here
we apply Theorem \ref{T6}.
\end{proof}

\begin{theorem}
\label{T4p} Under the assumptions of Theorem \ref{T8} we have%
\begin{eqnarray*}
\left| S(f;a,b)-P_{n}(f;a,b)\right| &\leq &\sqrt{\frac{47}{23040}-\frac{%
\sqrt{2}}{768}}\frac{(b-a)^{2}}{n^{2}}\sigma _{n}(f^{\prime }) \\
&\leq &\sqrt{\frac{47}{23040}-\frac{\sqrt{2}}{768}}\frac{(b-a)^{2}}{n\sqrt{n}%
}\omega _{n}(f^{\prime }),
\end{eqnarray*}%
where $S(f;a,b)$ and $P_{n}(f;a,b)$ are defined by (\ref{g0}) and (\ref{k01}%
), respectively and $\left\{ a=x_{0}<x_{1}<\cdots <x_{n}=b\right\} $ is a
uniform subdivision of $\left[ a,b\right] $, i.e. $x_{i}=a+ih,$ $h=(b-a)/n$, 
$i=0,1,...,n$.
\end{theorem}

\begin{proof}
The proof of this theorem is similar to the proof of Theorem \ref{T2p}. Here
we apply Theorem \ref{T8}.
\end{proof}

Nenad Ujevi\'{c}

Department of Mathematics

University of Split

Teslina 12/III

21000 Split

CROATIA

E-mail: ujevic@pmfst.hr


\begin{thebibliography}{99}
\bibitem{AT1} K. Atkinson, W. Han, Theoretical numerical analysis - A
functional analysis framework, Springer-Verlag, New York / Berlin /
Heidelberg, 2001.

\bibitem{C1} P. Cerone, Three points rules in numerical integration,
Nonlinear Anal.-Theory Methods Appl. 47 (4), (2001), 2341--2352.

\bibitem{CN1} D. Cruz-Uribe and C. J. Neugebauer, Sharp error bounds for the
trapezoidal rule and Simpson's rule, J. Inequal. Pure Appl. Math., 3(4),
Article 49, (2002), 1--22.

\bibitem{DAC1} S. S. Dragomir, R. P. Agarwal and P. Cerone, On Simpson's
inequality and applications, J. Inequal. Appl., 5 (2000), 533--579.

\bibitem{DCR1} S. S. Dragomir, P. Cerone and J. Roumeliotis, A new
generalization of Ostrowski's integral inequality for mappings whose
derivatives are bounded and applications in numerical integration and for
special means, Appl. Math. Lett., 13 (2000), 19--25.

\bibitem{DPW1} S. S. Dragomir, J. Pe\v{c}ari\'{c} and S Wang, The unified
treatment of trapezoid, Simpson and Ostrowski type inequalities for
monotonic mappings and applications, Math. Comput. Modelling, 31 (2000),
61--70.

\bibitem{GO1} A. Ghizzetti and A. Ossicini, Quadrature formulae, Birkh\"{a}%
user Verlag, Basel/Stuttgart, 1970.

\bibitem{KU1} A. R. Krommer, C. W. Ueberhuber, Computational integration,
SIAM, Philadelphia, 1998.

\bibitem{MPF1} D. S. Mitrinovi\'{c}, J. Pe\v{c}ari\'{c} and A. M. Fink,
Classical and new inequalities in analysis, Kluwer Acad. Publ.,
Dordrecht/Boston/London, 1993.

\bibitem{NS1} S. G. Nash and A. Sofer, Linear and Nonlinear Programing,
McGraw-Hill Book Co., New York/ Singapore, 1996.

\bibitem{PPUV1} C. E. M. Pearce, J. Pe\v{c}ari\'{c}, N. Ujevi\'{c} and S.
Varo\v{s}anec, Generalizations of some inequalities of Ostrowski-Gr\"{u}ss
type, Math. Inequal. Appl., 3(1), (2000), 25--34.

\bibitem{U5} N. Ujevi\'{c}, Inequalities of Ostrowski-Gr\"{u}ss type and
applications, Appl. Math., 29(4), (2002), 465--479.

\bibitem{U2} N. Ujevi\'{c}, Two sharp Ostrowski-like inequalities and
applications, (submitted).
\end{thebibliography}
\end{document}